\documentclass[12pt,a4paper]{article}

\usepackage[french]{babel}
\usepackage{amsmath}
\usepackage{amsfonts}
\usepackage{eufrak}
\usepackage{amssymb}
\usepackage{psfig}
\usepackage{epsf}
\usepackage{rotate}
\usepackage[all]{xy}
\usepackage{xy}
\xyoption{all}                 

\newcommand{\C}{\mathbb{C}}
\newcommand{\OO}{\mathcal{O}_{C}}
 \newcommand{\Z}{\mathbb{Z}}
\newtheorem{defi}{ D\'efinition}[section]
\newtheorem{theo}[defi]{ Th\'eor\`eme}
\newtheorem{prop}[defi]{ Proposition}

\newtheorem{rema}[defi]{ Remarque}
\newtheorem{lemm}[defi]{ Lemme}

\begin{document}
 \setcounter{page}{1}

 \title{Stabilit\'e des fibr\'es $\Lambda^{p}E_{L}$ et condition de Raynaud.}

\author{SCHNEIDER Olivier \thanks{Laboratoire J.-A. Dieudonn\'e
U.M.R. no 6621 du C.N.R.S.
 Universit\'e de Nice - Sophia Antipolis
 Parc Valrose
 06108 Nice Cedex 02
 France.
email : \texttt{oschneid@math.unice.fr}}}
\date{}

 \maketitle
\begin{abstract} 
Soit  $C$  une courbe lisse de genre $g \geq 2$ sur $\C$. Soit  $L$ un fibr\'e en droites sur $C$  engendr\'e par ses sections globales et $E_{L}$ le fibr\'e dual du noyau du morphisme d'\'evaluation  $e_{L}$. 
On \'etudie ici la relation entre la stabilit\'e et la v\'erification d'une condition $(R)$ introduite par Raynaud : on d\'emontre que lorsque $C$ est g\'en\'erale, $E_{L}$ est semi-stable ; on prouve ensuite que $E_{L}$ v\'erifie $(R)$ lorsque $\deg(L) \geq 2g$ ou bien lorsque $L$ est g\'en\'erique. Enfin on d\'emontre que pour tout $p$ dans $\{2,\dots, \mathrm{rg}(E_{L})-2\}$, si $\deg(L) \geq 2g+2$, $\Lambda^{p}E_{L}$ ne v\'erifie pas $(R)$.  

$$\mbox{\bf Abstract}$$
Let $C$ be a smooth curve of genus $g \geq 2$ on $\C$. Let $L$ be a line bundle on $C$ generated by its global sections and let  $E_{L}$ be the dual of the kernel of the evaluation map $e_{L}$. We are studying here the relation between the stability the fact that the bundle is verifying a condition $(R)$ introduced by Raynaud : we prove that $E_{L}$ is semi stable when $C$ is general. We also prove that $E_{L}$ is verifying $(R)$ when
$\deg(L) \geq 2g$ or when  $L$ is generic. Finally we prove that for each $p$ in $\{2,\dots, \mathrm{rg}(E_{L})-2\}$, if $\deg(L) \geq 2g+2$ then $\Lambda^{p}E_{L}$ is not verifying $(R)$. 

\end{abstract}

 \section{Introduction}
  Soit  $C$  une courbe lisse de genre $g \geq 2$ sur $\C$. Soit $J$ la Jacobienne de $C$. Dans ~\cite{R}, M.Raynaud
  introduit la condition suivante : un fibr\'e vectoriel $E$ sur $C$ v\'erifie $(\star)$
  si  il existe un ouvert non vide $U$ de $J$ tel que pour tout $L$ dans $U$,  $$\min(h^{0}(E \otimes L),h^{1}(E \otimes L))=0.$$
  On remarque que
  si  $E$ est un fibr\'e vectoriel sur $C$ de pente $g-1$, v\'erifiant $(\star)$, alors $E$ est semi-stable. On introduit ici
  une condition $(R)$ (qui implique $(\star)$) de fa\c{c}on \`a pouvoir \'etendre cette propri\'et\'e \`a tous les fibr\'es de pente enti\`ere.
 Les fibr\'es construits par M.Raynaud fournissent comme pour $(\star)$,
  des exemples de fibr\'es stables ne v\'erifiant pas $(R)$.
  Soit maintenant $L$ un fibr\'e en droites sur $C$  engendr\'e par ses sections globales, soit $M_{L}$ le noyau du morphisme d'\'evaluation
  $$e_{L}:H^{0}(C,L) \otimes \OO \longrightarrow L.$$
 On note $E_{L}:=M_{L}^{\ast}$. Lorsque $\deg(L) \geq 2g+1$, $E_{L}$ est stable (voir ~\cite{E}) et si $\deg(L)=2g$, alors $E_{L}$ est stable si et seulement si $L$ est tr\`es ample (voir ~\cite{B}). On \'etablira ici que pour $\deg(L) < 2g$, $E_{L}$ est semi-stable lorsque $C$ est g\'en\'erale. On \'evoquera en effet des cas de courbes sur lesquelles ces fibr\'es ne le sont pas. 
On verra ensuite que comme pour la semi-stabilit\'e, lorsque $L$ est g\'en\'erique sur une courbe g\'en\'erale,  $E_{L}$ v\'erifie $(R)$. Enfin on conclura en \'evoquant le cas des fibr\'es $\Lambda^{p}E_{L}$ lorsque $\deg(L) \geq 2g+2$ : ces fibr\'es sont semi-stables et on montrera qu'ils ne v\'erifient pas $(R)$.

 \section{Condition de Raynaud et semi-stabilit\'e}
 Dans tout ce qui suit $C$ est une courbe lisse de genre $g \geq 2$ sur $\C$. Pour tout  entier $n$, on note $J^{n}$ la vari\'et\'e
 qui param\'etrise les fibr\'es en droites de degr\'e $n$ sur $C$ ($J$ si $n=0$).
 Soit $E$ un fibr\'e vectoriel sur $C$
 de rang $r$ et de degr\'e $d$. On dit que $E$ v\'erifie la condition $(R)$ si :
 \begin{center} \hspace{0.5cm} $\forall n \in \mathbb{Z}$, pour $L$ g\'en\'erique dans $J^{n}$, $h^{0}(E \otimes L)$ ou $h^{1}(E \otimes L)$ est nul. \hfill{$(R)$}
 \end{center}
   On a les propri\'et\'es suivantes :
 \begin{itemize}
  \item[$\bullet$]  tout fibr\'e en droites v\'erifie $(R)$.
     \item[$\bullet$]  si $E$ v\'erifie $(R)$, alors :
                 \begin{itemize} \item         pour tout fibr\'e en droites $L$ sur $C$, $E\otimes L$ v\'erifie $(R)$.
                        \item   Le dual $E^{\ast}$ de $E$, v\'erifie $(R)$.    \end{itemize}
  \end{itemize}

 Si $E$ est non nul, on d\'efinit sa pente par :
 $$\mu(E):=\frac{d}{r}.$$
 Rappelons qu'un fibr\'e vectoriel $E$ sur $C$ est dit \textbf{stable} (resp. \textbf{semi-stable}) si tout sous-fibr\'e propre de
 $E$ a une pente strictement inf\'erieure (resp. inf\'erieure) \`a $\mu(E)$.
 On a des propri\'et\'es analogues a celles \'evoqu\'ees pr\'ec\'edemment pour la condition $(R)$ : tout fibr\'e en droites
 est stable ; si $E$ v\'erifie l'une ou l'autre de ces propri\'et\'es de stabilit\'e, il en sera de m\^eme pour $E^{\ast}$
 et pour $E\otimes L$, avec $L$ un fibr\'e en droites quelconque sur $C$. \\
En fait, pour qu'un fibr\'e $E$ v\'erifie $(R)$ il faut et il suffit de v\'erifier les deux conditions suivantes :
 \begin{enumerate}
  \item [1)]$h^{1}(E \otimes L)=0$ pour $L$ fibr\'e en droites sur $C$ g\'en\'erique de degr\'e  $g-1-[\mu(E)]$.
 \item[2)]$h^{0}(E \otimes L)=0$ pour $L$ fibr\'e en droites sur $C$ g\'en\'erique de degr\'e  $g-1-\ulcorner \mu(E) \urcorner $.
 \end{enumerate}
  En effet, pour tout diviseur $D$ sur $C$ de degr\'e positif, pour tout fibr\'e en droites $M$ sur $C$, on a la suite exacte
$$0 \longrightarrow  E \otimes M  \longrightarrow  E \otimes M(D)   \longrightarrow E \otimes \mathcal{O}_{D}  \longrightarrow 0. $$
En \'ecrivant la suite exacte longue d'homologie, on obtient : $$h^{1}(C,E \otimes M(D)) \leq h^{1}(C,E \otimes M) $$
et   $$ h^{0}(C,E \otimes M) \leq h^{0}(C,E \otimes M(D)).$$
\begin{itemize}
\item Si $n \geq  g-1-[\mu(E)]$ alors tout fibr\'e en droites g\'en\'erique de degr\'e $n$ s'\'ecrit $M(D)$ avec $M$ un fibr\'e en droites
g\'en\'erique de degr\'e $g-1-[\mu(E)]$ et $D$ un diviseur de degr\'e positif.
Comme $$h^{1}(C,E \otimes M(D)) \leq h^{1}(C,E \otimes M),$$   la condition $1)$ implique que $h^{1}(E \otimes L)=0$
g\'en\'eriquement lorsque $L$ est fibr\'e en droites sur $C$ de degr\'e  $n \geq  g-1-[\mu(E)]$.
 \item  Si $n \leq  g-1-\ulcorner \mu(E) \urcorner$ alors tout fibr\'e en droites g\'en\'erique de degr\'e $g-1-\ulcorner \mu(E) \urcorner$ s'\'ecrit $M(D)$ avec $M$ un fibr\'e en droites
 g\'en\'erique de degr\'e $n$ et $D$ un diviseur de degr\'e positif.
 Comme    $$ h^{0}(C,E \otimes M) \leq h^{0}(C,E \otimes M(D)),$$ la condition $2)$ implique que $h^{0}(E \otimes L)=0$
 g\'en\'eriquement lorsque $L$ est un fibr\'e en droites sur $C$ de degr\'e  $n \leq  g-1-\ulcorner \mu(E) \urcorner$.
 \end{itemize}
 \begin{prop}
 Soit $C$ une courbe lisse de genre $g$ sur $\C$. Si on a la suite exacte de fibr\'es vectoriels sur $C$ suivante  
$$0 \longrightarrow F \longrightarrow E \longrightarrow Q \longrightarrow 0,$$ alors  \label{sta}
$$E \ \mathrm{verifie} \ (R) \Longrightarrow \mu(F) \leq \ulcorner \mu(E) \urcorner \ \ et \ \ \mu(Q) \geq [\mu(E)]. $$
\end{prop}

\begin{description}
  \item [Preuve :]
  pour tout fibr\'e en droites $L$ on a
   $$h^{0}(C,F \otimes L)\leq h^{0}(C,E \otimes L).$$
     De plus si $\deg(L)=n$, on par Riemann-Roch
     $$h^{0}(C,E \otimes L)-h^{1}(C,E \otimes L)=r(E)(n+\mu(E)-(g-1)).$$
  De la m\^eme fa\c{c}on,
    $$h^{0}(C,F \otimes L)-h^{1}(C,F \otimes L)=r(F)(n+\mu(F)-(g-1)).$$
  Si $E$ v\'erifie $(R)$, alors pour $n=g-1- \ulcorner \mu(E) \urcorner$,
  $$h^{0}(C,F \otimes L)=h^{0}(C,E \otimes L)=0.$$
  Donc  $\chi( F \otimes L) \leq 0$, ce qui entra\^ine la premi\`ere in\'egalit\'e. \\
De la m\^eme fa\c{c}on, pour tout fibr\'e en droites $L$ on a 
$$h^{1}(C,Q \otimes L)\leq h^{1}(C,E \otimes L),$$ et $$h^{0}(C,Q \otimes L)-h^{1}(C,Q \otimes L)=r(Q)(n+\mu(Q)-(g-1)).$$
Alors pour $n=g-1-[ \mu(E)]$, 
$$h^{1}(C,Q \otimes L)=h^{1}(C,E \otimes L)=0.$$ 
Donc $\chi( Q \otimes L) \geq 0$, ceci entra\^ine la deuxi\`eme in\'egalit\'e.  \hfill{$\Box$}
\end{description}

 On voit  gr\^ace \`a ce r\'esultat que si $\mu(E)$ est entier, la condition $(R)$ entra\^ine la semi-stabilit\'e. Par contre la stabilit\'e n'entra\^ine pas la v\'erification de la condition $(R)$ (voir la construction de  Raynaud dans
  ~\cite{R}). 

\section{Stabilit\'e des fibr\'es $E_{L}$}

    Soit $L$ un fibr\'e en droites sur $C$ engendr\'e par ses sections globales ;
     soit  $M_{L}$ le  fibr\'e vectoriel de rang $h^{0}(C,L)-1$ sur $C$, noyau du morphisme d'\'evaluation :  \\
$$\hspace{2.4cm} \xymatrix{ 0 \ar[r]  &  M_{L}  \ar[r] & H^{0}(C,L)\otimes \OO \ar[r]^-{ev_{L}} & L \ar[r] &0 } \hspace{2.4cm} (1)$$
 On d\'efinit $E_{L}:=M_{L}^{\ast}$.
      L.Ein et R.Lazarsfeld d\'emontrent dans ~\cite{E}, que si $\deg(L) \geq 2g+1$, alors  $E_{L}$ est stable.
 Le cas $deg(L)=2g$ est trait\'e par A.Beauville dans ~\cite{B} :
 \begin{theo}(A.Beauville)
     Si $\deg(L)=2g$, alors $E_{L}$ est semi-stable et poss\`ede un diviseur th\^eta. De plus $E_{L}$ est stable si et seulement si $L$ est tr\`es ample.
 \end{theo}
 Traitons maintenant les cas o\`u $L$ est un fibr\'e en droites de degr\'e inf\'erieur \`a $2g$ ; on rappelle la propri\'et\'e suivante (voir ~\cite{L}) :
   \begin{lemm}(R.Lazarsfeld)
       Soit $C$ une courbe lisse de genre $g$ sur $\C$. Soit $L$ un fibr\'e en droites sur $C$, de degr\'e $g+d$ ($d \in \Z$) et engendr\'e par ses sections globales.
       Soient $x_{1},\dots,x_{d+h^{1}(L)-1}$, des points
       distincts sur $C$ tels que $$h^{1}(L(-\sum x_{i}))=h^{1}(L).$$ On a alors une suite exacte
         $$0 \longrightarrow \bigoplus_{i=1}^{d+h^{1}(L)-1} \OO(x_{i}) \longrightarrow E_{L} \longrightarrow L(-\sum_{i=1}^{d+h^{1}(L)-1} x_{i})
   \longrightarrow 0.$$  Ceci induit la suite exacte suivante pour tout $0 \leq p \leq d+h^{1}(L)-1$ :
  $$0 \longrightarrow \bigoplus_{i_{1}<\dots<i_{p}} \OO(x_{i_{1}}+\dots+x_{i_{p}}) \longrightarrow \Lambda^{p}E_{L}
  \longrightarrow  \bigoplus_{j_{1}<\dots<j_{d+h^{1}(L)-p}} L(-x_{j_{1}}-\dots-x_{j_{d+h^{1}(L)-p}})  \longrightarrow 0.$$  \label{c}
     \end{lemm}

   \begin{lemm}
    Soit $C$ une courbe lisse de genre $g$ sur $\C$. Soit $L$ un fibr\'e en droites sur $C$ engendr\'e par ses sections globales.
  Si $Q$ est un fibr\'e quotient propre de $E_{L}$ de rang $n$, alors   $$h^0(\det Q) \geq n+1 .$$   \label{n}
   \end{lemm}
  \begin{description}
 \item[Preuve du Lemme \ref{n} :]
 soit  $g+d$ le degr\'e de $L$. Si on note $h:=h^{1}(L)$, alors $r(E_{L})=d+h$.
     D'apr\`es le lemme \ref{c}, pour tous $x_{1},\dots,x_{d+h-1}$ sur $C$ tels que
   $h^{1}(L(-\sum x_{i}))=h$, on a
    $$0 \longrightarrow \bigoplus_{i_{1}<\dots<i_{n}} \OO(x_{i_{1}}+\dots+x_{i_{n}}) \longrightarrow \Lambda^{n}E_{L}
  \longrightarrow  \bigoplus_{j_{1}<\dots<j_{d+h-n}} L(-x_{j_{1}}-\dots-x_{j_{d+h-n}})  \longrightarrow 0.$$
   Comme $Q$ est un fibr\'e quotient de $E_{L}$ on a donc le diagramme suivant :
   $$\xymatrix{  0 \ar[d] \\
    \bigoplus_{i_{1}<\dots<i_{n}} \OO(x_{i_{1}}+\dots+x_{i_{n}})     \ar[d] \ar@{-->}[dr]| -{\varphi_{x_{1},\dots,x_{d+h-1}}} \\
    \Lambda^{n}E_{L} \ar[r] \ar[d] & \det Q \ar[r] & 0 \\
      \bigoplus_{j_{1}<\dots<j_{d+h-n}} L(-x_{j_{1}}-\dots-x_{j_{d+h-n}}) \ar[d]  \\
      0  }$$
  Dans tous les cas ceci impose $$h^0(\det Q) \geq n+1 .$$ En effet,
  \begin{itemize}
           \item[$\bullet$] ou bien il existe $x_{1},\dots,x_{d+h-1}$ tels que  $\varphi_{x_{1},\dots,x_{d+h-1}}$ est nulle ;
            alors il existe  $x_{i_{1}},\dots,x_{i_{d+h-n}}$ sur  $C$ tels que  $L(-x_{i_{1}}-\dots-x_{i_{d+h-n}})$ s'injecte dans $\det Q$. Or comme $$h=h^{1}(L) \leq h^{1}(L(-x_{i_{1}}-\dots-x_{i_{d+h-n}})) \leq h^{1}(L(-\sum x_{i}))=h,$$ on a par  Riemann-Roch :
                $$h^{0}(L(-x_{i_{1}}-\dots-x_{i_{d+h-n}}))=g+d-(d+h-n)-g+1+h=n+1 \leq  h^0(\det Q).$$
           \item[$\bullet$]Ou bien $\varphi_{x_{1},\dots,x_{d+h-1}}$ n'est jamais nulle : pour $x_{i_{1}},\dots,x_{i_{n}}$ g\'en\'eriques sur $C$,
                      $\OO(x_{i_{1}}+\dots+x_{i_{n}})$ s'injecte dans  $\det Q$ ; d'o\`u  $$h^0(\det Q) \geq n+1. $$     \hfill{$\Box$}
        \end{itemize}
  \end{description}

     \begin{prop}
Soit $C$ une courbe lisse  g\'en\'erale de genre $g\geq 3$ sur $\C$. Si $L$ est un fibr\'e en droites sur $C$ engendr\'e par ses sections globales, alors $E_{L}$ est semi-stable. \label{p}
 \end{prop}

  \begin{description}
 \item[Preuve :]
 soit  $g+d$ le degr\'e de $L$ et $h:=h^{1}(L)$. D\'emontrons la proposition par l'absurde : soit $Q$ un fibr\'e quotient
  propre de $E_{L}$ de rang $n<r(E_{L})=d+h$ et tel que $\mu(Q)<\mu(E_{L})$. Ceci impose
  $$\hspace{4.2cm} \deg(Q)< n \left( \frac{g+d}{d+h} \right)  \hspace{4.2cm} (2)$$
   D'apr\`es le lemme \ref{n}, $h^0(\det Q) \geq n+1$.
  Soit  $\rho$ le nombre de Brill-Noether pour les syst\`emes
   lin\'eaires de degr\'e $\deg(Q)$ et de dimension (projective) $n$ :
\begin{eqnarray} \rho & = & g-(n+1)(g-\deg(Q)+n)  \nonumber \\
                                      & = &-ng+(n+1)\deg(Q)-n(n+1). \nonumber \end{eqnarray}
 En utilisant $(2)$, on obtient :
   \begin{eqnarray} \rho& < & -ng+n(n+1)\left( \frac{g+d}{d+h}\right) -n(n+1) \nonumber \\
                                      & < & n\left( \frac{(n+1)(g-h)}{d+h}-g\right) \leq 0. \nonumber \end{eqnarray}
  Ceci contredit l'existence d'un tel fibr\'e $Q$ sur un courbe g\'en\'erale, ce qui prouve la semi-stabilit\'e. \hfill{$\Box$}
  \end{description}

  \begin{rema}
    \em  Le r\'esultat de la proposition n'est pas g\'en\'eralisable \`a toute courbe :
    soit  $C$ une courbe de genre $g \geq 3$,
 poss\'edant un syst\`eme lin\'eaire  $g_{d}^{1}$  (de degr\'e $d$ et de dimension projective $1$). Soit $M$ le fibr\'e
 en droites engendr\'e par ce syst\`eme lin\'eaire. Alors pour tout  fibr\'e en droites $L$ sur $C$  g\'en\'erique de degr\'e
 $g+d$, on a une injection $M \hookrightarrow L$, d'o\`u une surjection
 $$ E_{L} \longrightarrow E_{M} \longrightarrow 0.$$
 De plus
 $$\mu(E_{L})=\frac{g}{d}+1 \ \ \ et \ \ \ \mu(E_{M})=d.$$
 Alors si $d < \frac{g}{d}+1$, $E_{L}$ n'est pas semi-stable.   \label{d}

      \end{rema}

      \section{$E_{L}$ et la condition $(R)$}
    Soit  $L$ un fibr\'e en droites sur $C$ engendr\'e par ses sections globales. Soit $M$ un autre fibr\'e en droites sur $C$.
   Si on tensorise la suite exacte $(1)$ par $M$ et qu'on \'ecrit la suite
  exacte longue d'homologie, on obtient :  \\
$$\xymatrix{ 0 \ar[r] & H^{0}(C,M_{L}\otimes M) \ar[r] &
  H^{0}(C,L)\otimes H^{0}(C,M) \ar[r]^{\ \ \ \ \mu_{L,M}} & H^{0}(C,L\otimes M)  \ar[r] & \\
   \ar[r] & H^{1}(C,M_{L}\otimes M) \ar[r] & \dots \hspace{2cm} }$$   \\
  D'apr\`es ce qu'on a remarqu\'e pr\'ec\'edement, $E_{L}$ v\'erifiera la condition $(R)$
  si et seulement si les deux conditions suivantes sont v\'erifi\'ees :
  \begin{enumerate}
  \item [1)] $\mu_{L,M}$ est surjective pour $M$ fibr\'e en droites sur $C$ g\'en\'erique de degr\'e  $g-1-[\mu(M_{L})]$.
 \item[2)]$\mu_{L,M}$ est injective pour $M$ fibr\'e en droites sur $C$ g\'en\'erique de degr\'e  $g-1-\ulcorner \mu(M_{L}) \urcorner$.
 \end{enumerate}
   Remarquons tout d'abord que si $h^{0}(C,L) \leq 2$, alors $\mu_{L,M}$  est toujours de rang maximum ; en
   effet,
   \begin{itemize}
   \item ou bien  $ h^{0}(C,L)=1$  et $\mu_{L,M}$ est  trivialement toujours injective.
   \item ou bien $h^{0}(C,L)=2$
   et le "base-point free pencil trick" (voir ~\cite{ACGH} p152) nous donne le r\'esultat.
   \end{itemize}

  \begin{prop}
      Soit $C$ une courbe lisse de genre $g$ sur $\C$. Si $L$ est un fibr\'e en droites sur $C$ sans point base et de degr\'e superieur \`a $2g$, alors $E_{L}$ v\'erifie $(R)$.
\end{prop}
\begin{description}
 \item[Preuve :]
   si $\deg(L) \geq 2g$, alors $[\mu(M_{L})]=-2$. Pour que $E_{L}$ v\'erifie $(R)$ il faut donc que :
        \begin{enumerate}
  \item [1)] $\mu_{L,M}$ soit surjective pour $M$ fibr\'e en droites sur $C$ g\'en\'erique de degr\'e  $g+1$, et que
 \item[2)]$\mu_{L,M}$ soit injective pour $M$ fibr\'e en droites sur $C$ g\'en\'erique de degr\'e  $g$.
 \end{enumerate}
      La deuxi\`eme assertion est \'evidente ($h^{0}(C,M)=1$). Pour le premier point, comme $h^{0}(C,M)=2$ pour $M$ g\'en\'erique de degr\'e $g$, on  applique encore le "base-point free pencil trick" pour obtenir le r\'esultat.  \hfill{$\Box$}
\end{description}
Enfin comme sur une courbe g\'en\'erale, pour deux fibr\'es en droites g\'en\'eriques l'application $\mu_{L,M}$ est de rang maximum (voir ~\cite{Ba}), on a la proposition suivante : 
\begin{prop}
  Soit $C$ une courbe lisse g\'en\'erale de genre $g$ sur $\C$. Si $L$ est un fibr\'e en droites sur $C$ g\'en\'erique et engendr\'e par ses sections globales alors $E_{L}$ v\'erifie $(R)$.
\end{prop}

\begin{rema}  \em
   Pour les m\^emes raisons que pour la Proposition \ref{p}, ce r\'esultat n'est pas g\'en\'eralisable \`a toute courbe : en effet, reprenons l'exemple de la Remarque \ref{d} : $C$ est une courbe lisse de genre $g$ poss\'edant un $g_{d}^{1}$ avec $d+1 < \frac{g}{d}+1$. Soit $M$ le fibr\'e  en droites engendr\'e par ce syst\`eme lin\'eaire. Alors pour tout  fibr\'e en droites $L$ sur $C$  g\'en\'erique de degr\'e
 $g+d$, on a une surjection
 $$ E_{L} \longrightarrow E_{M} \longrightarrow 0,$$
et $\ulcorner \mu(E_{M}) \urcorner = \mu(E_{M})=d <  \mu(E_{L})=\frac{g}{d}+1$. De ce fait d'apr\`es la Proposition ~\ref{sta}, $E_{L}$ ne v\'erifie pas $(R)$. 
\end{rema}
 
\section{Le cas des fibr\'es $\Lambda^{p}E_{L}$ lorsque $\deg(L) \geq 2g+2$} 
Lorsque $\deg(L)\geq 2g+1$, $E_{L}$ est stable. De ce fait, pour tout $1<p<\mathrm{rg}(E_{L})$ $\Lambda^{p}E_{L}$ est polystable, c'est \`a dire somme directe de fibr\'es stables de m\^eme pente. 
Soit $\gamma:=\left[ \frac{g+1}{2} \right]$ ; dans ~\cite{P} M.Popa d\'emontre que lorsque $\deg(L) \geq g( \gamma +1)$, $\Lambda^{\gamma}E_{L}$ ne v\'erifie pas $(R)$. De ceci il d\'eduit le r\'esultat suivant (voir \cite{P}) :
 \begin{theo} (M.Popa)
          Pour tout $g \geq 2$, il existe un entier $\rho(g)$ tel que pour tout $r>\rho(g)$, il existe sur toute courbe de genre $g$
          un fibr\'e stable de rang $r$ ne v\'erifiant pas $(R)$.         \label{popa}

 \end{theo}
En s'inpirant de cela, on d\'emontre la Proposition suivante : 
\begin{prop}
 Soit $C$ une courbe lisse de genre $g \geq 2 $ sur $\C$. Si $L$ est un fibr\'e en droites sur $C$ de degr\'e superieur ou \'egal \`a $2g+2$, alors pour tout $p$ dans $\{2,\dots,\mathrm{rg}(E_{L})-2\}$, $\Lambda^{p}E_{L}$ ne  v\'erifie  pas $(R)$.
\end{prop}

\begin{description}
\item[Preuve : ]
soit $g+d$ le degr\'e de $L$, un calcul facile donne : 
$$\mu(\Lambda^p E_{L})=p+\frac{g}{d}p.$$
Remarquons tout d'abord que d'apr\`es le Lemme \ref{c}, pour $x_{1},\dots,x_{d-1}$ $d-1$ points g\'en\'eriques sur $C$, on a pour tout $p$ dans $\{2,\dots,d-2\}$, $$\OO(x_{1}+\dots+ x_{p}) \hookrightarrow \Lambda^{p}E_{L}.$$
Donc pour tous diviseurs effectifs $A_{q}$ et $B_{p}$ respectivement de degr\'es $q$ et $p$, $h^{0}(C,\Lambda^{p} E_{L}(A_{q}-B_{p})) \neq 0$. Etablissons maintenant gr\^ace \`a ceci les conditions sur $p$ pour que $\Lambda^{p}E_{L}$ ne v\'erifie pas $(R)$ : 
\begin{enumerate} 
\item[1)] Si $p \leq g$, tout fibr\'e en droites $M$ g\'en\'erique de degr\'e $g-2p$ est un fibr\'e engendr\'e par un diviseur de la forme $A_{g-p}-B_{p}$ avec $A_{g-p}$ et $B_{p}$ des diviseurs effectifs de degr\'es respectivement $g-p$ et $p$. D'apr\`es ce qui pr\'ec\`ede, $h^{0}(C,\Lambda^{p} E_{L} \otimes M) \neq 0$ et comme $$\chi(C,\Lambda^{p} E_{L}\otimes M) \leq 0 \Longleftrightarrow \mu(\Lambda^{p} E_{L}\otimes M) \leq g-1,$$ $\Lambda^{p}E_{L}$ ne v\'erifiera pas $(R)$ si  $p+\frac{g}{d}p+g-2p \leq g-1$, c'est \`a dire si 
 $p \geq \frac{d}{d-g}$.  
\item[2)] Si $p \geq g$ alors tout fibr\'e en droites $M$ g\'en\'erique de degr\'e $-p$ s'\'ecrit $\OO(-A_{p})$ avec $A_{p}$ diviseur effectif de degr\'e $p$. D'apr\`es ce qui pr\'ec\`ede, $h^{0}(C,\Lambda^{p} E_{L} \otimes M)\neq 0$ et pour les m\^emes raisons que dans le $1)$,  $\Lambda^{p}E_{L}$ ne v\'erifiera pas $(R)$ si 
$ p+\frac{g}{d}p-p \leq (g-1)$, c'est \`a dire si $p \leq d-\frac{d}{g}$.
\end{enumerate}
D'apr\`es les propri\'et\'es \'etablies dans le premier paragraphe, pour prouver la Proposition il nous suffit  de d\'emontrer que  $\Lambda^{p}E_{L}$ ou bien   $\Lambda^{d-p}E_{L}$ (qui est le dual de $\Lambda^{p}E_{L}$ \`a $\otimes L$ pr\`es) ne v\'erifie pas $(R)$ : 
\begin{enumerate}
\item[ i)] Si $d \leq 2g$ alors
\begin{itemize}
\item ou bien $\max(p,d-p) \geq g$, disons $p$ et alors d'apr\`es $2)$ $\Lambda^{p}E_{L}$ ne v\'erifie pas $(R)$ si $p \leq d-\frac{d}{g}$, ce qui est bien le cas puisque $d \leq 2g$.
\item ou bien $\max(p,d-p)  \leq g$ et  si on suppose que  $\Lambda^{p}E_{L}$ et    $\Lambda^{d-p}E_{L}$ v\'erifient $(R)$, d'apr\`es le $1)$, on doit avoir : $$p <d-\frac{d}{g} \ \ \mathrm{et} \ \  d-p<d-\frac{d}{g}.$$
Or ceci n'est pas possible lorsque $d \geq g+2$.
\end{itemize}
\item[ii)] Si $d \geq 2g$ alors 
\begin{itemize}
\item ou bien $\min(p,d-p) \leq g$ et d'apr\`es $1)$, on doit avoir $\min(p,d-p) \geq  \frac{d}{d-g}$  ce qui est automatique lorsque $d \geq 2g$.
\item  ou bien $p$ et $d-p$ sont sup\'erieurs ou \'egaux \`a $g$ et si on suppose que  $\Lambda^{p}E_{L}$ et    $\Lambda^{d-p}E_{L}$ v\'erifient $(R)$, d'apr\`es le $2)$, on doit avoir : $$p >d-\frac{d}{g} \ \mathrm{et} \ d-p>d-\frac{d}{g} .$$ Or cela implique $g\leq1$, \c{c}a n'est donc pas possible. \hfill{$\Box$}

\end{itemize}
  \end{enumerate}
\end{description}

 \section*{Remerciements}
Je remercie mon directeur de th\`ese Arnaud Beauville de m'avoir guid\'e dans ce travail de recherche.

\end{document}